\documentclass[a4paper,11pt]{amsart}
\usepackage{amssymb,amsfonts,amsxtra}
\usepackage[all]{xy}
\usepackage{fullpage}

\def\mathscr#1{\mathcal{#1}}

\newtheorem{theorem}{Theorem}[section]
\newtheorem{cor}[theorem]{Corollary}
\newtheorem{lem}[theorem]{Lemma}
\newtheorem{prop}[theorem]{Proposition}

\theoremstyle{definition}
\newtheorem{defi}[theorem]{Definition}

\newtheorem{rem}[theorem]{Remark}

\numberwithin{equation}{section}

\def\ds{\displaystyle}
\def\:{\colon}
\def\.{\cdot}

\def\<{\left\langle}
\def\>{\right\rangle}
\def\({\left(}
\def\){\right)}

\def\epsilon{\varepsilon}
\def\phi{\varphi}
\def\subset{\subseteq}

\def\leq{\leqslant}
\def\geq{\geqslant}

\def\bar#1{\overline{#1}}
\def\hat#1{\widehat{#1}}
\def\tilde#1{\widetilde{#1}}

\DeclareMathOperator*{\holim}{holim}

\def\Times_#1{{\ds\mathop{\times}_{#1}}}
\def\oTimes_#1{{\ds\mathop{\otimes}_{#1}}}
\def\Smash_#1{{\ds\mathop{\wedge}_{#1}}}

\def\k{\Bbbk}

\begin{document}
\title{Towers of $MU$-algebras and the generalized Hopkins-Miller
theorem}
\author{A.Lazarev }
\address{Mathematics Department, University of Bristol, Bristol, BS8 1TW,
England.} \email{a.lazarev@bristol.ac.uk}
\begin{abstract}
Our results are of three types. First we describe a general
procedure of adjoining polynomial variables to $A_\infty$-ring
spectra whose coefficient rings satisfy certain restrictions.A
host of examples of such spectra is provided by killing a regular
ideal in the coefficient ring of $MU$, the complex cobordism
spectrum. Second, we show that the algebraic procedure of
adjoining roots of unity carries over in the topological context
for such spectra. Third, we use the developed technology to
compute the homotopy types of spaces of strictly multiplicative
maps between suitable $K(n)$-localizations of such spectra. This
generalizes the famous Hopkins-Miller theorem and gives
strengthened versions of various splitting theorems.
\end{abstract}
 \maketitle
\begin{center}
\end{center}
\section{Introduction}
Our goals in this paper are three-fold.
First, we generalize and extend the methods of \cite{Laz} of constructing
 $A_\infty$ ring spectra by adjoining `polynomial variables'. More precisely, we work
with the categories of strict rings and modules constructed
\cite{EKMM} and adopt the terminology of the cited reference. In
particular we consider $S$-algebras and commutative $S$-algebras
rather than the equivalent notions of $A_\infty$ and $E_\infty$
ring spectra. So, suppose that $R$ is a commutative $S$-algebra
such that its coefficient ring $R_\ast$ has no elements of odd
degree and the graded ideal $I_\ast$ in $R_\ast$ is generated by a
regular sequence $u_1,u_2,\ldots$. Then it is known from work of
Strickland, \cite{Str} that the $R$-module $R/I$ obtained by
killing the ideal $I$ has a structure of an $R$-ring spectrum. We
assume that $R/I$ is actually an $R$-algebra. Then it turns out
that, informally speaking, all $R$-modules `between' $R/I$ and $R$
also possess structures of $R$-algebras. More precisely, denote by
$R/I(u_k^l)$ the $R$-module obtained from $R$ by killing the
sequence $(u_1,u_2,\ldots,u_{k-1},u_k^l,u_{k+1}, \ldots)$ in
$R_*$. Then we show that there are no obstructions for the
existence of strictly associative products on $R(u_k^l)$. The
natural `reduction maps' between these $R$-algebras also admit
strictly multiplicative liftings.

The strictly multiplicative products thus constructed are
typically not unique. In fact (unless there are some additional
assumptions such as the sparseness of the coefficient ring) they
are not unique even up to homotopy as shown by Strickland. There
is an infinite tower of obstructions to uniqueness in which the
first term is precisely Strickland's obstruction.

The basic example for our theory is provided by taking $R=MU$, the
complex cobordism $S$-algebra and $MU/I$ be the Eilenberg-MacLane
$S$-algebra $H\mathbb{Z}/p$. It follows by induction that all
$MU$-algebras obtained by killing any sequence of polynomial
generators and/or a prime $p$ have structures of $MU$-modules.

Second, we consider the question of adjoning roots of unity to an
$S$-algebra. This problem was also treated in the recent work
\cite{SVW} of Schwanzl, Vogt and Waldhausen. Their definition of a
topological extension has better formal properties than ours but
it applies in a far less general situation. We show that one can
adjoin roots of unity to such spectra as Morava $K$-theories
$K(n)$, Johnson-Wilson theories $E(n)$ and many other algebras
over the complex cobordism $S$-algebra $MU$.

Third, we address the problem of computing $S$-algebra maps
between a certain completion $\hat{E}(n)$ of $E(n)$ and an
$MU$-algebra $E$ which is assumed to be `strongly $K(n)$-complete'
 in some precise sense explained later on in the paper. Examples
of such $MU$-algebras include $\hat{E}(n)$, $K(n)$ and the
Artinian completion of the $v_n$-localization of the
Brown-Peterson spectum $BP$. It turns out that the space of
$S$-algebra maps between $\hat{E}(n)$ and a strongly
$K(n)$-complete $MU$-algebra is homotopically discrete with the
set of connected components being equal to the set of
multiplicative cohomology operations $\hat{E}(n)\rightarrow E$. We
call the results of this type the generalized Hopkins-Miller
theorem because the original Hopkins-Miller theorem
(cf.\cite{Rez}) asserts that the space of $S$-algebra self-maps of
the spectrum $E_n$ is homotopically equivalent to the (discrete)
space of multiplicative operations from $E_n$ to itself. Here
$E_n$ is a $2$-periodic version of the completed Johnson-Wilson
theory $\hat{E}(n)$ which came to be popularly known as the Morava
$E$-theory.

We chose to work with the $2(p^n-1)$-periodic theory $\hat{E}(n)$
rather than with $E_n$. The relation of $\hat{E}(n)$ to $E_n$ is
the same as the relation of the Adams summand  of $p$-completed
complex $K$-theory spectrum $\widehat{KU}$ to $\widehat{KU}$
itself. The advantage of $\hat{E}(n)$ is that it is smaller than
$E_n$, however it does not admit the action of the full Morava
stabilizer group.

One  consequence of our generalized Hopkins-Miller theorem is that $\hat{E}(n)$
admits a unique $S$-algebra structure, the result previously obtained in
\cite{Bak}.

Another consequence is that $\hat{E}(n)$
 splits off $E$ as an $S$-algebra
for a certain class of strongly $K(n)$-complete $MU$-algebras $E$.
Such splittings were
previously known to be multiplicative only up to homotopy.

Throughout the paper the symbol $\mathbb{F}_p$ will denote the
prime field with $p$ elements, $L/\mathbb{F}_p$ an arbitrary (but
fixed) finite separable extension of $\mathbb{F}_p$. Further
$W(L)$ and $W_n(L)$ denote the ring of Witt vectors and the ring
of Witt vectors of length $n$ respectively. For an $R$-algebra $A$
we denote by $A^e$ its enveloping $R$-algebra $A\wedge_RA^{op}$
where $A^{op}$ is the algebra $A$ with opposite multiplication.

{\bf Acknowledgement.} The author is indebted to the referee for a
host of very valuable suggestions and comments. Thanks are also
due to A. Baker for many useful discussions on the subject of the
paper.

\section{Adjoining polynomial variables to $S$-algebras}
In this section we assume that $R$ is a fixed commutative $S$-algebra
such that $R_\ast=\pi_\ast R$ is a graded commutative ring concentrated in even degrees.
 We will also assume without loss of generality that
$R$ is  $q$-cofibrant in the sense of \cite {EKMM}. All objects under consideration
will be $R$-modules or $R$-algebras and smash products and homotopy classes of maps will
be understood to be taken in the category of $R$-modules.

We begin by reminding the reader the notions of topological derivations and
topological singular extensions of $R$-algebras. A detailed account can be
found in \cite{Laz}.
  Let $A$ be an  $R$-algebra and $M$ an $A$-bimodule. Then the
$R$-module $A\vee M$ has the obvious structure of an $R$-algebra
(`square-zero extension' of $A$). Consider the set $[A,A\vee
M]_{R-alg/A}$ of homotopy classes of $R$-algebra maps from $A$ to
$A\vee M$ which commute with the projection onto $A$. Then there
exists an $A$-bimodule $\Omega_A$ and a natural in $M$ isomorphism
$$[A,A\vee M]_{R-alg/A}\cong [\Omega_A,M]_{A-bimod}$$
 where the right hand side denotes the homotopy classes of maps in the category
of $A$-bimodules.
\begin{defi}
The topological derivations  $R$-module of $A$ with values in $M$
is the function $R$-module
$F_{A\wedge A^{op}}(\Omega_A,M)$.
We denote it by ${\bf Der}(A,M)$ and its $i$th homotopy group
by $Der^{-i}(A,M)$.
\end{defi}
The $A$-bimodule $\Omega_A$ is constructed as the homotopy fibre
of the multiplication map $A\wedge A\rightarrow A$. There exists
the following homotopy fibre sequence of $R$-modules:
\begin{equation}\label{sdf}{\bf THH}(A,M)\rightarrow M \rightarrow
{\bf Der}(A,M)\end{equation} Here ${\bf THH}(A,M)$ is the
topological Hochschild cohomology spectrum of $A$ with values in
$M$, ${\bf THH}(A,M):=F_{A\wedge A^{op}}(A,M)$.

We will frequently use the notion of a primitive operation from
$A$ to $M$. Denote by $m:A\wedge A\rightarrow A$ the
multiplication map and by $m_l:A\wedge M\rightarrow M$ and
$m_r:M\wedge A\rightarrow M$ the left and right actions of $A$ in
$M$ respectively. Then a map $p:A\rightarrow M$ is called
primitive  if $p$ is a `derivation up to homotopy', i.e. the
following diagram is homotopy commutative: $$\xymatrix{A\wedge
A\ar[r]^m\ar[d]_{1\wedge p\vee p\wedge 1}&A\ar[d]^p\\ A\wedge
M\vee M\wedge A\ar[r]^-{m_l\vee m_r}&M}$$ There is a forgetful map
$l:Der^\ast(A,M)\rightarrow [A,M]^\ast$ defined as follows. For
any topological derivation $d:A\rightarrow A\vee M$ let $l(d)$ be
the composite map $A\rightarrow A\vee M\rightarrow M$ where the
last map is just the projection onto the wedge summand. Then it is
easy to see that the image of $l$ is contained in the set of
primitive operations from $A$ to $M$.

Suppose we are given a topological derivation $\xi:A\rightarrow
A\vee M$. Consider the following homotopy pullback diagram
$$\xymatrix{X\ar[r]\ar[d]&A\ar[d]\\ A\ar[r]^-{\xi}&A\vee M}$$ Here
the rightmost downward arrow is the canonical inclusion of a
retract. Then we have the following homotopy fibre sequence of
$R$-modules:
\begin{equation}\label{er}\Sigma^{-1}M\rightarrow X\rightarrow A\end{equation}
\begin{defi} The homotopy fibre sequence (\ref{er}) is called the topological
singular extension associated with the topological derivation
$\xi:A\rightarrow A\vee M$.
\end{defi}

For an element $x\in R_*$ we will denote by $R/x$ the
 cofibre of the map $R\stackrel{x}{\rightarrow}R$.
Let $I_\ast$ be a graded ideal generated by (possibly infinite)
regular sequence of homogeneous elements $(u_1,u_2,\ldots)\in
R_\ast$. We assume in addition that each $u_k$ is a nonzero
divisor in the ring $R_*$. (This assumption is satisfied in all
cases of interest). Then we can form the $R$-module $R/I$ as the
infinite smash product of $R/u_k$. By \cite{Str}, Proposition
$4.8$ there is a structure of an $R$-ring spectrum on $A$. Clearly
the coefficient ring of $R/I$ is isomorphic to $R_\ast/I_\ast$
where $R_\ast/I_\ast$ is understood to be the direct limit of
$R_\ast/(u_k,u_2,\ldots,u_k)$.

Let us denote the $R$-algebra $R/I$ by $A$. Our standing
assumption is that $A$ has a structure of an $R$-algebra (i.e.
strictly associative). This may seem a rather strong condition
but, as we see shortly such a situation is quite typical. In fact
P.Goerss proved in \cite{Goe} that any spectrum obtained by
killing a regular sequence in $MU$, the complex cobordism
$S$-algebra, has a structure of an $MU$-algebra.

The construction we are about to describe allows one to construct
new $R$-algebras by `adjoining' the indeterminates $u_k$ to the
$R$-algebra $A$. The basic idea is the same as in \cite{Laz} where
Morava $K$-theories at an odd prime were shown to possess
$MU$-algebra structures. However the arguments we use here are
considerably more general, in particular we make no assumption
that the prime $2$ is invertible.
 Let us introduce
the notation $A(u^l_k)_*$ for the $R_\ast$-algebra
$$\lim_{n\rightarrow \infty}R_\ast/(u_1,u_2,
\ldots,u_{k-1},u_k^l,u_{k+1},\ldots).$$
 For each $l$ reduction modulo
$u_k^{l}$ determines a map of $R_\ast$-algebras
$${A(u^{l+1}_k)_*}\rightarrow {A(u^l_k)_*}. $$ Now we can
formulate our main theorem in this section.
\begin{theorem} \label{qas}
For each $l$ there exist $R$-algebras $A(u^l_k)$ with coefficient
rings ${A(u^l_k)_*}$ and $R$-algebra maps
$$R_{l,k}:A(u^{l+1}_k)\rightarrow A(u^l_k)$$ which give the
reductions mod $u_k^{l}$ on the level of coefficient rings.
\end{theorem}
{\bf Proof}. Set $d:=|u_k|$. Suppose by induction that the
$R$-algebras $A(u^l_k)$ with the required properties were
constructed for $l\leq i$. We will show that there exists an
appropriate Bokstein operation from  $A(u^i_k)$ to
$\Sigma^{di+1}A$ which allow us to build the next stage.

Consider the cofibre sequence
$$\Sigma^{dl}R\stackrel{u_k^l}{\rightarrow}R\stackrel{\rho_{l,k}}{\rightarrow}
R/u_k^l\stackrel{\beta_{l,k}}\rightarrow\Sigma^{dl+1}R.$$
According to \cite{Str} the $R$-module $R/u_k^l$ admits an
associative product $\phi:R/u_k^{l}\wedge R/u_k^{l}\rightarrow
R/u_k^{l}$ and for any other product $\phi^{\prime}$ there exists
a unique element $u\in \pi_{2dl+2}R/u_k^{l}$ for which
$\phi^\prime=\phi+u\circ(\beta_{l,k}\wedge\beta_{l,k})$.

\begin{lem} There exists a map of $R$-ring spectra
$ r_{i,k}:R/u_k^{i+1}\rightarrow R/u_k^{i}$
 realizing
the reduction map  on coefficient rings. Moreover in the cofibre
sequence
\begin{equation}\label{bk1}
R/u_k^{i+1}\stackrel{r_{i,k}}\rightarrow
R/u_k^{i}\stackrel{\bar{\beta}_{i,k}} {\rightarrow}
\Sigma^{di+1}R/u_k\end{equation} the second map
$\bar{\beta}_{i,k}: R/u_k^i\rightarrow \Sigma^{di+1}R/u_k$ is a
primitive operation.
\end{lem}
{\bf Proof.} Use induction on $i$ (in fact we only need the
inductive assumption in order that $R/u_k$ be an algebra spectrum
over $R/u_k^i$.) Consider the following diagram of $R$-modules:
$$\xymatrix{\Sigma^{d(i+1)}R\ar[r]^{u_k}\ar[d]_{u_k^{i+1}}&\Sigma^{di}R
\ar[r]\ar[d]_{u_k^{i}}&\Sigma^{di}R/u_k\ar[d]\\R\ar@{=}[r]\ar[d]&R\ar[r]\ar[d]&\mbox{pt}
\ar[d]
\\R/u^{i+1}_k\ar[r]^{r_{i,k}}&R/u^i_k\ar[d]\ar[r]&\Sigma^{di+1}R/u_k\ar@{=}[d]\\
&\Sigma^{di+1}R\ar[r]&\Sigma^{di+1}R/u_k}$$ Here the rows and
columns are homotopy cofibre sequences of $R$-modules. The map
$r_{i,k}:R/u_k^{i+1}\rightarrow R/u_k^{i}$ is determined uniquely
by the requirement that the diagram commute in the homotopy
category of $R$-modules.

We now show that it is possible to choose a product on
$R/u_k^{i+1}$ so that $r_{i,k}$ becomes an $R$-ring spectrum map.
First take any associative product $\phi:R/u_k^{i+1}\wedge
R/u_k^{i+1}\rightarrow R/u_k^{i+1}$ which exists by \cite{Str},
Proposition $3.1$. As in \cite{Str}, Proposition $3.15$ there is
an obstruction $d(\phi)\in\pi_{2d(i+1)+2}(R/u_k^{i})$ for the map
$r_{i,k}$ to be homotopy multiplicative. If the product $\phi$ is
changed the obstruction changes according to the formula
$$d(\phi+u\circ(\beta_{i,k}\wedge\beta_{i,k}))=d(\phi)+r_{i,k\ast}u.$$
Since the map $r_{i,k\ast}:R_\ast/u_k^{i+1}\rightarrow
R_\ast/u_k^{i}$ is just the reduction map it is surjective and we
conclude that there exists a product on $R/u_k^{i+1}$ for which
the obstruction vanishes.

Further the map $\bar{\beta}_{i,k}: R/u_k^i\rightarrow
\Sigma^{di+1}R/u_k$ coincides with the composition
$$R/u_k^i\stackrel{\beta_{i,k}}{\rightarrow} \Sigma^{di+1}R
\stackrel{\rho_{i,k}}{\rightarrow}\Sigma^{di+1}R/u_k^i\rightarrow
\Sigma^{di+1}R/u_k.$$ Here the last map is (a suspension of) the
canonical $R$-ring map $R/u_k^i\rightarrow R/u_k$ which exists by
the inductive assumption. The composition
$\rho_{i,k}\circ\beta_{i,k}$ is the Bokstein operation which is
primitive by \cite{Str}, Proposition $3.14$. It follows that
$\bar{\beta}_{i,k}$ is also primitive and our lemma is proved.

Now taking the smash product of (\ref{bk1}) with $R/u_1\wedge
R/u_2\wedge\ldots
 \wedge R/u_{k-1}\wedge R/u_{k+1}\wedge\ldots$ we get the following homotopy cofibre sequence
of $R$-modules: $$A(u^{i+1}_k)\stackrel{R_{i,k}}\longrightarrow
A(u^i_k) \stackrel{Q_{i,k}}\longrightarrow \Sigma^{di+1}A.$$ It
follows that $Q_{i,k}$ is a primitive operation and $A(u^{i+1}_k)$
has an $R$-ring spectrum structure such that $R_{i,k}$ is an
$R$-ring map. (Notice our abuse of notations here in using the
symbol $R_{i,k}$ even though $A(u^{i+1}_k)$ is not yet proved to
be an $R$-algebra).

Next we will describe the set of all primitive operations
$A(u^{i}_k)\rightarrow \Sigma^\ast A$. Consider the cofibre
sequence $$\Sigma^{|u_l|}R\stackrel{u_l}\rightarrow
R\stackrel{\rho_l}{\rightarrow}R/u_l
\stackrel{\beta_l}{\rightarrow}\Sigma^{|u_l|+1}R.$$ For $l\neq k$
we introduce the operation $$Q_l:A(u^{i}_k)= R/u_1\wedge
R/u_2\wedge \ldots\wedge R/u_{k-1}\wedge R/u^i_k\wedge
R/u_{k+1}\wedge \ldots\rightarrow \Sigma^{|u_l|+1}A$$ obtained by
smashing $\rho_l\circ\beta_l:R/u_l\rightarrow
\Sigma^{|u_l|+1}R/u_l$ with the identity map on the remaining
smash factors. The operations $Q_l$ are primitive, again by
\cite{Str}, Proposition 3.14. The next result shows that $Q_l$ and
$Q_{i,k}$ are essentially all primitive operations from $A(u^i_k)$
into (suspensions of) $A$.
\begin{lem} Any primitive operation $A(u^i_k)\rightarrow \Sigma^\ast A$
can be written uniquely as an infinite sum
$a_kQ_{i,k}+\Sigma_{i\neq k}a_iQ_i$ for $a_i\in A_\ast$.
\end{lem}
{\bf Proof.} We will only give a sketch since the arguments of
Strickland (\cite{Str}, Proposition 4.17) carry over almost
verbatim. Let $\eta:R\rightarrow A(u^i_k)$ be the unit map and
$$J_\ast=(u_1,u_2,\ldots,u_{k-1},u^i_k,u_{k+1},\ldots)$$ be the
kernel of the map $\eta_\ast:R_\ast\rightarrow {A(u^i_k)}_\ast$.
Given a primitive operation $Q:A(u^i_k)\rightarrow\Sigma^sA$
define the function $p(Q):J_\ast\rightarrow A(u^i_k)_\ast$ as
follows. For any $x\in J_\ast$ we have a cofibre sequence
$$\Sigma^{|x|}R\stackrel{x}{\rightarrow}R\stackrel{\rho_x}{\rightarrow}R/x
\stackrel{\beta_x}{\rightarrow} \Sigma^{|x|+1}R.$$ Then there is a
unique map $f_x:R/x\rightarrow A(u^i_k)$ such that $f_x\circ
\rho_x=\eta$. Further there is a unique map
$y:\Sigma^{|x|+1}R\rightarrow \Sigma^sA$ such that $Q\circ
f_x=y\circ \beta_x$. We define $p(Q)(x):=y\in
\pi_{|x|+1-s}A_\ast$.

It follows that the function $p$ actually embeds the set of
primitive operations  $A(u^i_k)\rightarrow \Sigma^*A$ into
$Hom^*_{R_\ast}(J/J^2, A_\ast)$. Further it is straightforward to
check that $$p(Q_s)(u_t)=\delta_{st} ~\mbox{for}~ s,t\neq k,~
p(Q_s)(u_k^i)=0;$$ $$p(Q_{i,k})(u_j)=0 ~\mbox{for}~ j\neq k,~
p(Q_{i,k})(u_k^i)=1.$$ Therefore the elements $Q_{i,k}, Q_s, s\neq
k$ form a basis in $Hom_{R_*}(J_\ast/J^2_\ast,A_*)$ dual to the
basis $(u_1,u_2,\ldots,u_{k-1},u_k^i,u_{k+1}\ldots)$ in
$J_\ast/J^2_\ast$ and our lemma is proved.

Now we come to the crucial part of the proof. We are going to show
that the Bokstein operation $Q_{i,k}$ can be improved to a
topological derivation from which it would follow that
$A(u^{i+1}_k)$ is an $R$-algebra and
$R_{i,k}:A(u^{i+1}_k)\rightarrow A(u^i_k)$ lifts to an $R$-algebra
map.

The following lemma (interesting in its own right) is preparatory
for computing the topological Hochschild cohomology of
$A(u^{i}_k)$ with coefficients in $A$.
 \begin{lem}\label{opp} For any regular ideal $J_\ast=(x_1,x_2,\ldots)$ in $R_\ast$
and any product on $R/J$ there is  a multiplicative isomorphism
$$\pi_\ast(R/J)^e\cong \Lambda_{R_\ast/J_\ast}(\tau_1,
\tau_2,\ldots)$$ where $|\tau_i|=|x_i|+1$.
\end{lem}
{\bf Proof.} We have the universal coefficients spectral sequence
$$E^2_{**}=Tor^{R_*}_{**}(R_*/J_*,(R_*/J_*)^{op})\Longrightarrow
\pi_{*}(R/J)^e.$$ This is a standard right half-plane spectral
sequence of homological type with $E^\infty$ term associated to an
increasing filtration $$R_*/J_*=F_0\subset F_1\subset\ldots\subset
\pi_*(R/J).$$ This spectral sequence is multiplicative as shown in
\cite{BaL}, Lemma 1.3. Therefore using the standard Koszul
resolution for the $R_*$-module $R_*/J_*$ we see that
$E^2_{**}=\Lambda_{R_*}(\tau_1,\tau_2,\ldots).$ The differentials
of this spectral sequence vanish on the generators $\tau_i$ and we
conclude that $E^2_{**}=E^\infty_{**}$.

We will now show that there are no nontrivial multiplicative
extensions in $E^*_{\ast\ast}$ This is not immediate since it is
not a spectral sequence of {\it commutative} algebras. Take a
representative of the generator $\tau_i$ in $\pi_*(R/J)^e$. The
indeterminacy in choosing this representative is an odd degree
element in $R_*/J_*$. Since the ring $R_*/J_*$ is even this
representative is in fact determined canonically and we will still
denote it by $\tau_i$.  Since the associated graded ring to
$\pi_*(R/J)^e$ is graded commutative we conclude that $\tau_i^2$
belongs to the filtration component $F_1$. Since $F_1/F_2$ is
spanned by $\tau_i$'s which are odd-dimensional it follows that in
fact $\tau_i\in F_0=R_*/J_*$.

Next consider the map of $R$-ring spectra $f:(R/J)^e\rightarrow
F_R(R/J,R/J)$ which is induced by the structure of an
$R/J$-bimodule spectrum on $R/J$. Then $f$ induces a map of
$R_\ast/J_\ast$-algebras $f_\ast:\pi_\ast (R/J)^e\rightarrow
\pi_\ast F_R(R/J,R/J)$. By \cite{Str}, Proposition $4.15$ the ring
$\pi_\ast F_R(R/J,R/J)$ is a completed exterior algebra over
$R_\ast/J_\ast$. Therefore $f_*(\tau_i^2)=0$. Since $\tau^2\in
R_*/J_*$ and the map $f$ is $R_*/J_*$-linear it implies that
$\tau_i^2=0$. Similarly all graded commutators of the elements
$\tau_i$ vanish in the ring $\pi_*(R/J)^e$. This finishes the
proof of Lemma \ref{opp}.
\begin{rem} The $R_\ast/J_\ast$-algebra $\pi_\ast(R/J\wedge_RR/J)$ need not
be an exterior algebra in general. For example take $R/J=MU/2$,
the reduction of $MU$ modulo $2$. Then it can be shown using the
methods of \cite{Nas} that $\pi_*MU/2\wedge_{MU}MU/2$ is an
$MU/2_\ast$-algebra on one generator in degree $1$ whose square is
equal to $x_1\in MU/2_\ast$. We will discuss this and related
phenomena elsewhere.
\end{rem}
\begin{lem}\label{tor}
There is the following isomorphism of graded $R_\ast$-modules:
$$grTHH^\ast(A(u^i_k), A)=
A_\ast[\tilde{u}_1,\tilde{u}_2,\ldots,\tilde{u}_k,\ldots];$$
$$grDer^{\ast-1}(A(u^i_k), A)=
A_\ast[\tilde{u}_1,\tilde{u}_2,\ldots,\tilde{u}_k,\ldots]/(A_\ast).$$
Here $gr(?)$ denotes the associated graded module. Moreover the
image of  the forgetful map $$Der^\ast(A(u^i_k), A)\rightarrow
[A(u^i_k), A]^\ast$$ is the set of all primitive operations from
$A(u^i_k)$ to $\Sigma^\ast A$.
\end{lem}
{\bf Proof.} Denote by $J_\ast$ the ideal in $R$ generated by
$(u_1,u_2,\ldots,u_{k-1},u_k^i,u_{k+1},\ldots)$. Then we have the
following spectral sequence $$Ext^{\ast\ast}_{\pi_\ast
(A(u^i_k))^e} (A(u^i_k)_\ast,A_\ast) =Ext^{\ast\ast}_{\pi_\ast
(R/J)^e} (R/J_\ast,A_\ast) \Rightarrow THH^{\ast}(A(u^i_k), A).$$
By Lemma \ref{opp} $\pi_\ast
(R/J)^e=\Lambda_{R/J_\ast}(\tau_1,\tau_2\ldots)$. Therefore
$$Ext^{\ast\ast}_{\pi_\ast (R/J)^e}
(R/J_\ast,A)=A_\ast[\tilde{u}_1,\tilde{u}_2,\ldots,\tilde{u}_k,\ldots]$$
where $|\tilde{u}_l|=-|u_l|+2$ for $l\neq k$ and
$|\tilde{u}_k|=-di+2$. Since our spectral sequence is even it
collapses and we obtain the desired isomorphism
$$grTHH^\ast(A(u^i_k), A)=
A_\ast[\tilde{u}_1,\tilde{u}_2,\ldots,\tilde{u}_k,\ldots].$$ The
isomorphism involving $Der^\ast(A(u^i_k), A)$ is obtained
similarly.

Now consider the spectral sequence
$$Ext^{\ast\ast}_{R_*}(A(u^i_k)_\ast, A_\ast)=
Hom^*(\Lambda(z_1,z_2,\ldots), A_\ast)\Rightarrow [A(u^i_k),
A]^\ast.$$ Here $|z_l|=|u_l|+1$ for $l\neq k$ and $|z_k|=di+1$.
This spectral sequence is easily seen to collapse and denoting  by
$\tilde{z_i}$ the elements in the dual basis in
$\Lambda(z_1,z_2,\ldots)$ we identify its $E_2=E_\infty$-term with
$\Lambda_{A_\ast}(\tilde{z}_1,\tilde{z}_2,\ldots)$. (Of course
this is an identification only as $A_\ast$-modules, not as rings.)
It is clear that the elements $z_l$ correspond to the primitive
operations $Q_l:A(u^i_k)\rightarrow \Sigma^{|u_l|+1}A$ for $l\neq
k$ whilst $z_k$ corresponds to $Q_{i,k}:A(u^i_k)\rightarrow
\Sigma^{|d|i+1}A$.

Furthermore the forgetful map $l$ (operating on the level of
$E_2$-terms) sends the element $\tilde{u}_l\in
grDer^\ast(A(u^i_k), A)$ to $\tilde{z}_l$. In other words all
primitive operations $A(u^i_k)\rightarrow \Sigma^\ast A$ are
covered by $l$ up to higher filtration terms. Since the image of
$l:Der^\ast(A(u^i_k), A)\rightarrow [A(u^i_k), A]^\ast$ is
contained in the subspace of the primitive operations no higher
filtration terms are present and we conclude that all Bokstein
operations $Q_l$ and $Q_{i,k}$ are in the image of $l$. With this
Lemma \ref{tor} is proved.
\begin{rem}Part of the Lemma \ref{tor} could be reformulated more
canonically as follows: the associated graded to the filtered ring
$THH^*(A(u_k^i),A)$ is isomorphic to the symmetric algebra over
$A_*$ generated by the $A_*$-module $J/IJ$ (with degrees shifted
by one).
\end{rem}

 So we proved that there exists a
topological derivation $$\tilde{Q}_{i,k}:A(u^i_k)\rightarrow
A(u^i_k)\vee \Sigma^{di+1}A$$ such that its composition with the
projection onto the wedge summand $$
A(u^i_k)\vee\Sigma^{di+1}A\rightarrow \Sigma^{di+1}A$$ is the
Bokstein operation $Q_{i,k}$. Associated with $\tilde{Q}_{i,k}$ is
a topological singular extension $$\Sigma^{di}A\rightarrow
X\rightarrow A(u^i_k)$$ such that $X$ is weakly equivalent to
$A(u^{i+1}_k)$. In other words the $R$-module $A(u^{i+1}_k)$
admits a structure of an $R$-algebra so that the reduction map
$A(u^{i+1}_k)\rightarrow A(u^i_k)$ is an $R$-algebra map. The
inductive step is completed and Theorem \ref{qas} is proved.

The singular extension of $R$-algebras
$$A(u^{i+1}_k)\stackrel{R_{i,k}}\longrightarrow A(u^i_k)
\stackrel{Q_{i,k}}\longrightarrow \Sigma^{di+1}A$$ is of
fundamental importance to us. We will call the $R$-algebra
$A(u^{i+1}_k)$ an elementary extension of $A(u^i_k)$. We will
denote by $A(u_k^\infty)$ the $R$-algebra
$\holim_{i\rightarrow\infty}A(u^i_k)$. In the case of adjoining
several, or perhaps, infinitely many indeterminates
$u_{k_1},u_{k_2},\ldots$ we will write
$A(u_{k_1}^\infty,u_{k_2}^\infty,\ldots)$ for the $R$-algebra
$$\holim_{l\rightarrow\infty}A(u_{k_1}^\infty)(u_{k_2}^\infty)
\ldots (u_{k_l}^\infty).$$

Now let $R$= $MU$, the complex cobordism spectrum. It is
well-known that $MU$ is a commutative $S$-algebra and
$MU_\ast=\mathbb{Z}[x_1,x_2,\ldots]$, the polynomial algebra on
infinitely many generators in even degrees.
\begin{cor}\label{klo} Let $J_\ast$ be a regular ideal in $MU_\ast=\mathbb{Z}[x_1,x_2,\ldots]$
 generated by any  subsequence of the regular sequence
$p^{\nu_0},x_{1}^{\nu_1},x_{2}^{\nu_2},\ldots$
where $\nu_l\geq 1$ for all $l$. Then $MU/J$ admits a structure of an $MU$-algebra.

\end{cor}
{\bf Proof.} First assume that the element $p^{\nu_0}$ does not
belong to $J_\ast$. Denote by $I_\ast$ the ideal in $MU_\ast$
generated by all polynomial generators $(x_1,x_2,\ldots)$. Then
$MU/I$ is the integral Eilenberg-MacLane spectrum $H\mathbb{Z}$
which possesses a canonical structure of an $MU$-algebra (even a
commutative $MU$-algebra).

Taking successive elementary extensions corresponding to elements
$x_{i_{k}}^{\nu_{i_k}}$ with $\nu_{i_k}>1$ we construct an
$MU$-algebra $$\widetilde{MU/J} =
H\mathbb{Z}(x_{i_1}^{\nu_{i_1}},x_{i_2}^{\nu_{i_2}}\ldots).$$
Similarly taking elementary extensions corresponding to those
polynomial generators $\{u_1,u_2,\ldots\}$ whose powers are not in
$J_\ast$ and passing to the (inverse) limit we construct the
$MU$-algebra $\widetilde{MU/J}(u_1^\infty,u_2^\infty,\ldots)$
whose underlying $MU$-module is $MU/J$.

If the element $p^{\nu_0}$ does belong to $J_\ast$ the only
difference is that we start our induction with the
Eilenberg-MacLane $MU$-algebra $H\mathbb{Z}/p^{\nu_0}$ instead of
$H\mathbb{Z}$. With this Corollary \ref{klo} is proved.

\begin{rem}Notice that our method also shows that the $MU$-algebra structures
on $MU/J$ for different $J$ are compatible in the sense that various reduction
maps given by killing elements $x_{i}^{\nu_i}$ are actually $MU$-algebra maps.
\end{rem}
\begin{rem} The statement of Corollary \ref{klo} remains true with
$MU$ replaced with $MU_{(p)}$, the localization of $MU$ at any
prime $p$. The proof is the same verbatim. \end{rem}
 Recall that
the Brown-Peterson spectrum $BP$ is obtained from $MU_{(p)}$ by
killing all polynomial generators in $MU_{(p)*}$ except for
$v_i=x_{2(p^i-1)}$ provided that we use Hazewinkel generators for
$MU_{(p)*}=\mathbb{Z}[x_1,x_2\ldots,]$. Therefore
$BP_\ast=\mathbb{Z}_{(p)}[v_1,v_2,\ldots]$. It follows that $BP$
possesses an $MU$-algebra structure. We define $$BP\langle
n\rangle=BP/(v_i,i>n)=H\mathbb{Z}_{(p)}(v_1^\infty,v_2^\infty,\ldots,v_n^\infty)$$
$$P(n)=BP/(v_i,i<n)=H\mathbb{F}_p(v_n^\infty,v_{n+1}^\infty,\ldots)$$
$$B(n)=v_n^{-1}P(n)=v_n^{-1}[H\mathbb{F}_p(v_n^\infty,v_{n+1}^\infty,\ldots)]$$
$$k(n)=BP/(v_i,i\neq n)=H\mathbb{F}_p(v_n^\infty)$$
$$K(n)=v_n^{-1}BP/(v_i,i\neq
n)=v_{n}^{-1}[H\mathbb{F}_p(v_n^\infty)]$$
$$E(n)=v_n^{-1}BP/(v_i,i>
n)=v_{n}^{-1}[H\mathbb{Z}_{(p)}(v_1^\infty,v_2^\infty,\ldots,v_n^\infty)]$$
Note that inverting an element in the coefficient ring of an
$MU_{(p)}$-ring spectrum is an instance of Bousfield localization
in the homotopy category of $MU_{(p)}$-modules. Further Bousfield
localization preserves algebra structures we conclude that all
spectra listed above admit structures of $MU_{(p)}$-algebras (and
therefore $MU$-algebras).

We conclude this section with a few remarks about the commutativity of the products
on the $MU$-algebras considered. First in the case of the odd prime $p$
$BP$ and other spectra derived from it have coefficient rings concentrated in
degrees congruent to $0$ mod $4$ and therefore all products are unique up to
homotopy and automatically commutative. If $p$=2 then Strickland shows that
$BP$ still admits a structure of a commutative $MU$-ring spectrum, but this
does not follow directly from our construction. We conjecture
that  in the context of Corollary \ref{klo} one can choose
an $MU$-algebra structure on $MU/J$ compatible
with any given structure of an $MU$-ring spectrum on $MU/J$.

The situation with $E(n)$ is different. The spectrum $E(n)$ is
Landweber exact regardless of the prime $p$. Therefore it
possesses a unique homotopy associative and commutative
multiplication compatible with its structure (in the classical
sense) of an $MU$-algebra spectrum. Our results show that this
multiplication can be improved to an $MU$-algebra structure which
in turn gives an $A_\infty$-structure. On the other hand
Strickland shows in \cite{Str} that at the prime $2$ there is no
homotopy commutative multiplication on $E(n)$ in the category of
$MU$-modules (to get a homotopy commutative product one has to
choose a different set of generators). Therefore, quite curiously,
$E(n)$ is a homotopy commutative $A_\infty$-ring spectrum and
simultaneously a  $MU$-algebra that is not commutative even up to
homotopy.

\section{Separable extensions of $R$-algebras}
We keep our convention that $R$ is a commutative $S$-algebra with coefficient
ring $R_\ast$ concentrated in even degrees. Let $A$ be an $R$-algebra with coefficient
ring $A_\ast$ and $A_\ast\subset\overline{A}_\ast$ a ring extension in the usual
algebraic
sense. In this section we consider the problem of finding an $R$-algebra $\overline{A}$
together with an $R$-algebra map $A\rightarrow \overline{A}$ wich realizes the given
extension $A_\ast\subset\overline{A}_\ast$ in homotopy. If such an algebra $
\overline{A}$
exists we say that the algebraic extension $A_\ast\subset\overline{A}_\ast$ admits a
topological lifting.

We now describe a general framework in which it is possible to
prove that topological liftings exist. Suppose that
$R_\ast\rightarrow \overline{R}_\ast$ is a separable extension of
$R$. That means that $R_\ast$ is a subring of a graded commutative
ring $\overline{R}_\ast$ and $\overline{R}_\ast$ is a separable
algebra over $R_\ast$, i.e. $\overline{R}_\ast$ is a projective
module over $\overline{R}_\ast\otimes_{R_\ast}\overline{R}_\ast$.
In addition we assume that $\overline{R}_\ast$ is a projective
$R_\ast$-module. Then for any ideal $I_\ast$ in $R_\ast$ the
$R_\ast$-algebra
$\overline{R_\ast/I}_\ast:=R_\ast/I_\ast\otimes_{R_\ast}\overline{R}_\ast$
is a separable $R_\ast/I_\ast$-algebra so
$R_\ast/I_\ast\rightarrow \overline{R_\ast/I}_\ast$ is also a
separable extension.

Further assume that the ideal $I_\ast$ is generated by a regular
sequence $(u_1,u_2,\ldots)$ (possibly infinite) of nonzero
divisors, the $R$-module $A=R/I$ is supplied with a structure of
an $R$-algebra and the algebraic extension
$R_\ast/I_\ast\rightarrow \overline{R_\ast/I}_\ast$ admits a
topological lifting. In other words there exists an $R$-algebra
$\overline{A}$ and an $R$-algebra map $A\rightarrow  \overline{A}$
realizing the given separable extension on the level of
coefficient rings.

Recall that in the previous section we constructed an $R$-algebra
structure on the $R$-module $$A(u_k^l):=R/u_1\wedge
R/u_2\ldots\wedge R/u_{k-1}\wedge R/u_k^{l}\wedge
R/u_{k+1}\wedge\ldots.$$ Then we have the following
\begin{theorem} \label{nh}
There exist  $R$-algebras $\overline{A(u_k^l)}$ realizing in
homotopy the $R_\ast$-module $\overline{{A(u_k^l)}}_\ast:=
{A(u_k^l)}_\ast\otimes_{R_\ast}\overline{R}$ and $R$-algebra maps
$\overline{R_{l,k}}:\overline{A(u_k^{l+1})}
\rightarrow\overline{A(u_k^l)}.$ which give the reduction maps on
the level of coefficient rings. Moreover the algebraic extension
of rings ${A(u_k^l)}_\ast\rightarrow \overline{{A(u_k^l)}}_\ast$
admits a topological lifting so that the following diagram of
$R$-algebras is homotopy commutative: $$\xymatrix{
A(u_k^{l+1})\ar[r]^{R_{l,k}}\ar[d] & A(u_k^l)\ar[d]\\
\overline{A(u_k^{l+1})}\ar[r]^{R_{l,k}}& \overline{A(u_k^l)}}$$
\end{theorem}
Before embarking on the proof of the theorem we will formulate and
prove a useful lemma which will be called `flat base change'. It
will be used several times in this paper. Though the lemma is
really a standard exercise in homological algebra we have been
unable to find a proper reference and will therefore give a
complete proof.
\begin{lem} Let $k$ be a graded commutative ring and $B,C$ be
graded $k$-algebras such that $C$ is flat as a $k$-module. Let $M$
and be a graded left $C$-module which is flat as a $k$-module and
$N$ be a left graded ${C\otimes_kB}$-module. Then there is the
following `flat base change' isomorphism: $$Ext^{**}_C(M,N)\cong
Ext^{**}_{C\otimes_kB}(M\otimes_kB,N).$$\end{lem} {\bf Proof}. We
have the following change of rings spectral sequence corresponding
to the algebra map $id\otimes 1:C\rightarrow {C\otimes_kB}$ (cf.
\cite{CE}, Ch. XVI, 6):
\begin{equation}\label{hyper}Ext^{**}_{C\otimes_kB}(Tor_{**}^C(C\otimes_kB,M),N)\Rightarrow
Ext_C^{**}(M,N)\end{equation} In fact, this spectral sequence was
deduced in the cited reference for ungraded rings and modules, but
it does not affect the arguments. Since $M$ and $C$ are flat
$k$-modules we obtain $$Tor_{**}^C(C\otimes_kB,M)\cong
Tor_{**}^C(B\otimes_kC,M)\cong Tor_{**}^k(B,Tor^{**}_C(C,M)$$ $$
\cong B\otimes_CC\otimes_kB\cong B\otimes_kM\cong M\otimes_kB.$$
It follows the the spectral sequence (\ref{hyper}) collapses and
our lemma is proved.
\begin{cor}\label{hoh}Let $k\hookrightarrow \overline{k}$ be a separable
extension of a graded commutative ring $k$. Let $C$ be a graded
ring and $M,N$ be graded left $\overline{k}\otimes C$-modules.
Then there is the following isomorphism:
$$Ext_{C\otimes\overline{k}\otimes_k\overline{k}}(M\otimes
\overline{k} ,N)\cong Ext_C(M,N).$$\end{cor}{\bf Proof}. Since the
$\overline{k}\otimes_k\overline{k}$-module $\overline{k}$ is
projective it is flat. Moreover we obviously have
$\overline{k}\otimes_{\overline{k}\otimes_k\overline{k}}\overline{k}\cong
\overline{k}$. Therefore by Lemma \ref{hyper}

$$Ext_{C\otimes\overline{k}\otimes_k\overline{k}}(\overline{k}\otimes
M,N)\cong Ext_{C\otimes\overline{k}}
(\overline{k}\otimes_{\overline{k}\otimes_k\overline{k}}\otimes
\overline{k}\otimes M,N)$$ $$ \cong
Ext_{C\otimes\overline{k}}(M\otimes \overline{k},N)\cong
Ext_C(M,N)$$

{\bf Proof} of Theorem \ref{nh}. Proceeding by induction suppose
that the $R$-algebras $\overline{A(u_k^l)}$ with required
properties were constructed for $l\leq i$.
 Consider the spectral sequence
$$Tor_{\ast\ast}^{R_\ast}
({\overline{A(u_k^i)}}_\ast,{\overline{A(u_k^i)}}^{op}_\ast) \cong
\overline{R}_\ast\otimes_{R_\ast}\overline{R}_\ast
\otimes_{R_\ast} Tor_{\ast\ast}^{R_\ast}
({A(u_k^i)}_\ast,{A(u_k^i)}^{op}_\ast)$$
$$\cong\overline{R}_\ast\otimes_{R_\ast}\overline{R}_\ast
\otimes_{R_\ast} \Lambda_{{A(u_k^i)}_\ast}(\tau_1,\tau_2,\ldots)
\cong \overline{R}_\ast\otimes_{R_\ast}
\overline{A(u_k^i)}_\ast\otimes \Lambda(\tau_1,\tau_2,\ldots)
\Rightarrow \pi_\ast\overline{A(u_k^i)}^e.$$ Here the first
isomorphism is an instance of the isomorphism
$$Tor_*^R(P\otimes_RM,Q\otimes_RN)\cong
P\otimes_RQ\otimes_RTor_*^R(M,N)$$ for projective $R$-modules $M$
and $N$. The second isomorphism follows from the standard
calculation with the Koszul complex and the third isomorphism is
obvious.

 Since by Lemma \ref{opp}
$\pi_\ast
A(u_k^i)^e=\Lambda_{A(u_k^i)_\ast}(\tau_1,\tau_2,\ldots)$ we see
that the exterior generators $\tau_i$ are permanent cycles and it
follows that our spectral sequence collapses multiplicatively so
the ring $\pi_\ast\overline{A(u_k^i)}^e$ is isomorphic to the
exterior algebra on $\tau_1,\tau_2,\ldots$ with coefficients in
$\overline{R}_\ast\otimes_{R_\ast} \overline{A(u_k^i)}_\ast$.

Next notice that the $R$-algebra map $A(u^i_k)\rightarrow
\overline{A(u_k^i)}$ induces a map of $R$-modules
\begin{equation}\label{ret}THH(A(u^i_k),A(u^i_k))\rightarrow
THH(A(u^i_k),\overline{A(u^i_k)})\end{equation} We claim that
(\ref{ret}) is a weak equivalence. To see this consider the
spectral sequences \begin{equation}\label{ret1}
Ext^{**}_{A(u^i_k)^e_*}(A(u^i_k)_*,\overline{A(u^i_k)_*})\Rightarrow
THH^*(A(u^i_k),\overline{A(u^i_k)})\end{equation} and
\begin{equation}\label{ret2}
Ext^{**}_{\overline{A(u^i_k)^e_*}}(\overline{A(u^i_k)_*},
\overline{A(u^i_k)_*})\Rightarrow
THH^*(\overline{A(u^i_k)},\overline{A(u^i_k)})\end{equation} The
map (\ref{ret}) determines a map between spectral sequences
(\ref{ret1}) and (\ref{ret2}). From Corollary \ref{hoh} we deduce
that these spectral sequences are isomorphic and therefore
(\ref{ret}) is a weak equivalence.

Finally we conclude that the map ${\bf
Der}(\overline{A(u^i_k)},\overline{A}) \rightarrow{\bf
Der}(A(u^i_k),\overline{A})$ induced by the $R$-algebra map
$A(u^i_k)\rightarrow \overline{A(u^i_k)}$ is weak equivalence.

Consider the following homotopy commutative diagram of
$R$-algebras:
\begin{equation}\label{dia}\xymatrix{\overline{A(u^i_k)}\ar[d]^\xi&
A(u^i_k)\ar[l]\ar[d]&A(u^i_k)\ar@{=}[l]\ar[d]^{Q_{i,k}}\\
\overline{A}\vee
\Sigma^{di+1}\overline{A}&A\vee\Sigma^{di+1}\overline{A}
\ar[l]&A\vee\Sigma^{di+1}{A}\ar[l] }
\end{equation}
Here the vertical arrows are topological derivations. The right
downward arrow is the Bokstein operation $Q_{i,k}$ constructed in
the previous section, the middle arrow is the only one that makes
the right square commute. Since derivations of
$\overline{A(u^i_k)}$ with values in $\overline{A}$ are in
one-to-one correspondence with derivations of $A(u^i_k)$ with
values in $\overline{A}$ the left downward arrow $\xi$ making the
left square commute exists and is unique.

Now the derivation $\xi$ gives rise to a topological singular
extension $$\Sigma^{di}\overline{A}\rightarrow ?\rightarrow
\overline{A(u^i_k)}$$ which on the level of coefficient rings
reduces to the algebraic singular extension
$$\Sigma^{di}\overline{A}_\ast\rightarrow
\overline{A(u^{i+1}_k)}_\ast \rightarrow
\overline{A(u^i_k)}_\ast$$ Therefore we can denote $?$ by
$\overline{A(u^{i+1}_k)}$ and because of the diagram (\ref{dia})
we have a map of topological singular extensions
$$\xymatrix{\Sigma^{di}{A}\ar[d]\ar[r]&A(u^{i+1}_k)\ar[d]\ar[r]&
A(u^i_k)\ar[d]\\
\Sigma^{di}\overline{A}\ar[r]&\overline{A(u^{i+1}_k)}
\ar[r]&\overline{A(u^i_k)}}$$ With this the inductive step is
completed and Theorem \ref{nh} is proved.

We now show that Theorem \ref{nh} gives a way of adjoining roots
of unity to various $MU$-algebras. Let
$\mathbb{F}_p\hookrightarrow L$ be a finite separable extension of
the field $\mathbb{F}_p$ (typically obtained by adjoining roots of
an irreducible factor of the cyclotomic polynomial $x^{p^n-1}-1$
for some $n$). Let $R=\widehat{MU}_p$, the $p$-completion of the
complex cobordism spectrum $MU$ and $\overline{R}_\ast=W(L)\otimes
MU_\ast$ where $W(L)$ is the ring of Witt vectors of $L$. Since
$L$ is a separable extension of $\mathbb{F}_p$ the algebra $W(L)$
is separable over the $p$-adic integers $\hat{\mathbb{Z}}_p$ and
it follows that $W(L)\otimes MU_\ast$ is a separable algebra over
$\widehat{MU}_{p*}=\hat{\mathbb{Z}}_p\otimes MU_\ast$.

\begin{cor}\label{adj}
Let $J_\ast$ be the ideal in
$\widehat{MU}_{p\ast}=\hat{\mathbb{Z}}_p[x_1,x_2,\ldots]$
generated by $(p, x_{i_1},x_{i_2},\ldots)$ where $x_{i_k}$ are
polynomial generators in $\widehat{MU}_{p\ast}$. Then for any
finite separable extension $\mathbb{F}_p\hookrightarrow L$ there
exist $\widehat{MU}_p$-algebras $MU/J$, $MUL/J$, $\widehat{MUL/J}$
and $\widehat{MU/J}$ together with the commutative diagram in the
homotopy category of $\widehat{MU}_p$-algebras
\begin{equation}\label{qa}\xymatrix{\widehat{MU/J}\ar[r]\ar[d]&\widehat{MUL/J}\ar[d]\\
{MU/J}\ar[r]&MUL/J}\end{equation} which reduces in homotopy to the
following diagram of $MU_\ast$-algebras:
$$\xymatrix{\hat{\mathbb{Z}}_p\otimes
MU_\ast/(x_{i_1},x_{i_2}\ldots)\ar[r]\ar[d]_{\mbox{mod p
}}&W(L)\otimes MU_\ast/(x_{i_1},x_{i_2}\ldots)\ar[d]^{\mbox{mod p
}}\\ MU_\ast/J_\ast\ar[r]& L\otimes MU_\ast/J_\ast}$$
\end{cor}
{\bf Proof.} Denote by $u_1,u_2,\ldots$ the collection of
polynomial generators of $\widehat{MU}_{p\ast}$ which are not in
$J_\ast$. There is an isomorphism of rings $MU_\ast/J_\ast\cong
\widehat{MU}_{p\ast}/J_*\cong \mathbb{F}_p[u_1,u_2,\ldots]$. Let
$I$ be the maximal ideal $(p,x_1,x_2,\ldots)$ in
$\widehat{MU}_{p\ast}=\hat{\mathbb{Z}}_p[x_1,x_2,\ldots]$. Then
$\widehat{MU}_p/I$ is the Eilenberg-MacLane spectrum
$H\mathbb{F}_p$ and the canonical map $\widehat{MU}_p\rightarrow
MU/I$ is a commutative $S$-algebra map. Therefore $MU/I$ is an
$\widehat{MU}_p$-algebra (even a commutative
$\widehat{MU}_p$-algebra). Further the extension
$\mathbb{F}_p\hookrightarrow L$ determines a map of commutative
$S$-algebras $H\mathbb{F}_p\rightarrow HL$ which is also a map of
commutative $\widehat{MU}_p$-algebras. Using Theorem \ref{nh} we
can adjoin the variable $u_i$ to the topological extension
$H\mathbb{F}_p\rightarrow HL$ and obtain the tower of topological
extensions $$\{H\mathbb{F}_p(u^i_k)\rightarrow HL(u^i_k)\}.$$
 Passing to the limit we get the topological extension
$H\mathbb{F}_p(u_k^\infty)\rightarrow HL(u_k^\infty)$.

This procedure could be repeated so we can adjoin another
polynomial generator to the extension
$H\mathbb{F}_p(u_i^\infty)\rightarrow HL(u_k^\infty)$
 (or any number of polynomial
generators). In this way we construct the extension
$$H\mathbb{F}_p(u_1^\infty,u_2^\infty,\ldots)\rightarrow
HL(u_1^\infty,u_2^\infty,\ldots)$$ which is the same as the
extension $MU/J\rightarrow MUL/J$ for the ideal $J_\ast=(p,
x_{i_1},x_{i_2},\ldots)$ in $\widehat{MU}_{p\ast}$. At this point
we adjoin the indeterminate corresponding to the prime $p$. We
then get the following diagram of $\widehat{MU}_p$-algebras:
$$\xymatrix{MU/J\ar[d]&MU/J(p^2)\ar[d]\ar[l]&\ldots\ar[l]&MU/J(p^n)\ar[d]\ar[l]&
\ldots\ar[l]\ar[l]\\
MUL/J&MUL/J(p^2)\ar[l]&\ldots\ar[l]&MUL/J(p^n)\ar[l]&\ldots\ar[l]}$$
where the $\widehat{MU}_p$-algebras $MU/J(p^n)$ and $MUL/J(p^n)$
realize the $\widehat{MU}_{p\ast}$-algebras
$W_n(\mathbb{F}_p)\otimes MU_\ast/(x_{i_1},x_{i_2},\ldots)$ and
$W_n(L)\otimes MU_\ast/(x_{i_1},x_{i_2},\ldots)$ respectively.

Finally taking  the (homotopy) inverse limit we get an
$\widehat{MU}_p$-algebra map
$$MU/J(p^\infty)=\widehat{MU/J}\rightarrow
\widehat{MUL/J}=MUL/J(p^\infty)$$ which realizes in homotopy the
extension of rings $\hat{\mathbb{Z}}_p\otimes MU/J_\ast\rightarrow
W(L)\otimes MU/J_\ast.$ Corollary \ref{adj} is then proved.
\begin{rem} It is not hard to prove that the diagram (\ref{qa}) is equivariant
with respect to the Galois group $Gal(L/\mathbb{F}_p)$ in a
suitable sense but we save this observation for future work.
\end{rem}

\section{Localized towers of $MU$-algebras}
In this section we investigate the behaviour of the towers of $MU$-algebras
constructed in the previous section under Bousfield localization. We start with
an almost trivial
\begin{theorem}\label{wed}
Let \begin{equation}\label{sin}
I\rightarrow A\rightarrow B
\end{equation} be a singular extension of $R$-algebras
for a commutative $S$-algebra $R$. Then for any $R$-module $M$ the cofibre
sequence \begin{equation}\label{si}
I_M\rightarrow A_M\rightarrow B_M\end{equation} is also a singular extension
of $R$-algebras (here $?_M$ denotes Bousfield localization with respect to $M$)
\end{theorem}
{\bf Proof.} The singular extension (\ref{sin}) is associated with
a certain topological derivation $\xi: B\rightarrow B\vee \Sigma
M$ so that there is a homotopy pullback square of $R$-algebras
$$\xymatrix{ A\ar[d]\ar[r]&B\ar[d]\\ B\ar[r]^-\xi& B\vee \Sigma I
}$$ where the vertical map $B\rightarrow B\vee \Sigma M$ is the
canonical inclusion of the wedge summand. Localizing this  diagram
with respect to $M$ we get the homotopy pullback diagram
$$\xymatrix{A_M\ar[d]\ar[r]&B_M\ar[d]\\ B_M\ar[r]^-{\xi_M}&B_M\vee
\Sigma I_M }$$ Therefore (\ref{si}) is a singular extension
associated with the topological derivation
$B_M\stackrel{\xi_M}{\rightarrow} B_M\vee \Sigma I_M$ and the
theorem is proved.

Now let $I_\ast$ be a regular ideal $(u_{1}, u_{2},\ldots)$ in
$R_\ast$ and $A=R/I$ be an $R$-algebra. Consider the $R$-algebra
$A(u_1^\infty, u_{2}^\infty,\ldots,
u_{k}^\infty)_{u_1^{-1}[A(u_1^\infty]})$, the Bousfield
localization of the $R$-algebra $A(u_1^\infty,
u_{2}^\infty,\ldots, u_{k}^\infty)$ with respect to
$u_{1}^{-1}[A(u_1^\infty)]$. Then we have
\begin{prop}\label{oj}\emph{(i)} There is the following weak equivalence of
$R$-algebras: $$A(u_1^\infty, u_{2}^\infty,\ldots,
u_{k}^\infty)_{u_1^{-1}[A(u_1^\infty)]}\simeq
(u_1^{-1}[A(u_1^\infty)])(u_2^\infty,u_3^\infty,\ldots,u_k^\infty).$$
 \emph{(ii)} If $R$ is connective and the elements $u_1,u_2,\ldots, u_k$
have positive degrees then the coefficient ring of the $R$-algebra
$A(u_1^\infty, u_{2}^\infty,\ldots,
u_{k}^\infty)_{u_1^{-1}[A(u_1^\infty)]}$ is the $R_*$-algebra
$A_\ast[u_{1}^{\pm 1}][[u_{2}, u_{3},\ldots, u_{k}]]$.
\end{prop}
{\bf Proof.} For (i) consider the following tower of $R$-algebras:
$$A(u_1^\infty)\leftarrow A(u_1^\infty,u_2^2)\leftarrow\ldots
\leftarrow A(u_1^\infty,u_2^l)\leftarrow \ldots$$
Bousfield-localizing it with respect to $u_1^{-1}A(u_1^\infty)$ we
get the tower $$u_1^{-1}[A(u_1^\infty)]\leftarrow
u_1^{-1}[A(u_1^\infty,u_2^2)]\leftarrow\ldots \leftarrow
u_1^{-1}[A(u_1^\infty,u_2^l)]\leftarrow \ldots.$$ Notice that the
first term of the localized tower as well as its successive
subquotients (isomorphic to  suspensions of
$u_1^{-1}[A(u_1^\infty)]$) are $u_1^{-1}[A(u_1^\infty)]$-local and
therefore so is its homotopy inverse limit. Clearly the canonical
map $$\holim_{l\rightarrow\infty} A(u_1^\infty,u_2^l)\rightarrow
\holim_{l\rightarrow\infty}
u_1^{-1}A(u_1^\infty,u_2^l)=(u_1^{-1}[A(u_1^\infty)])(u_2^\infty)$$
is $u_1^{-1}[A(u_1^\infty)]$-equivalence and we conclude that
\begin{equation}\label{fu}(u_1^{-1}[A(u_1^\infty)])(u_2^\infty)\cong
A(u_1^\infty,u_2^\infty)_{u_1^{-1}[A(u_1^\infty)]}\end{equation}
We then proceed by adjoining $u_{i_3}$. There is a tower of
$R$-algebras $$A(u_1^\infty)(u_2^\infty)\leftarrow
A(u_1^\infty)(u_2^\infty,u_3^2) \leftarrow \ldots\leftarrow
A(u_1^\infty)(u_2^\infty,u_3^l)\leftarrow\ldots$$ where the
successive stages
$$\Sigma^{|u_{3}|l}A(u_1^\infty)(u_2^\infty)\rightarrow
A(u_1^\infty)(u_2^\infty,u_3^{l+1})\rightarrow
A(u_1^\infty)(u_2^\infty,u_3^l)$$ are singular extensions of
$R$-algebras. Localizing this tower with respect to
$u_1^{-1}[A(u_1^\infty)]$ and using (\ref{fu}) and Theorem
\ref{wed} we get the tower
$$(u_1^{-1}[A(u_1^\infty)])(u_2^\infty)\leftarrow
(u_1^{-1}[A(u_1^\infty)])(u_2^\infty,u_3^2) \leftarrow
\ldots\leftarrow
(u_1^{-1}[A(u_1^\infty)])(u_2^\infty,u_3^l)\leftarrow\ldots$$whose
homotopy limit $$\holim_{l\rightarrow \infty}
(u_1^{-1}[A(u_1^\infty)])(u_2^\infty,u_3^l)=
(u_1^{-1}[A(u_1^\infty)])(u_2^\infty,u_3^\infty)$$ is clearly
weakly equivalent  to the localization of
$A(u_1^\infty)(u_2^\infty,u_3^\infty)$ with respect to
$u_1^{-1}A(u_1^\infty)$.  Repeating this process for
$u_{4},u_{5},\ldots,u_{k}$ we obtain the desired isomorphism
$$A(u_1^\infty, u_{2}^\infty,\ldots,
u_{k}^\infty)_{u_1^{-1}[A(u_1^\infty)]}\simeq
(u_1^{-1}[A(u_1^\infty)])(u_2^\infty,u_3^\infty,\ldots,u_k^\infty).$$
For (ii) notice that since $u_1$ has positive degree the ring
$A(u_1^\infty)_*$ is a polynomial ring over $A$ on one variable
$u_1$. Furthermore $u_1^{-1}[A(u_1^\infty)]_*=A_*[u_{1}^{\pm 1}]$.
Denote by $B$ the ring
$$(u_1^{-1}[A(u_1^\infty)])(u_2^\infty,u_3^\infty,\ldots,u_k^\infty)_*\cong
A(u_1^\infty, u_{2}^\infty,\ldots,
u_{k}^\infty)_{u_1^{-1}[A(u_1^\infty)]*},$$ Then $B$ is by
definition the completion of the ring
$u_1^{-1}R_*/(u_{k+1},u_{k+2},\ldots)$ at the ideal
$(u_2,u_3,\ldots,u_{k})$. Since this ideal is regular the
associated graded ring is polynomial: $$grB=A_\ast[u_{1}^{\pm
1}][u_{2}, u_{3},\ldots, u_{k}].$$ Since the elements $u_l,
l=1,\ldots,k$ have positive degrees the $R_*$-algebra
$A(u_1^\infty, u_{2}^\infty,\ldots, u_{k}^\infty)_{*}$ is in fact
an $A_*$-algebra and therefore, so is $B$. It follows that $B$ is
isomorphic to the completion of its associated graded ring and our
proposition is proved.

Our main application of the developed localization techniques is
to the Johnson-Wilson theories $E(n)$. Recall that the
$MU$-algebra spectrum  $E(n)$ has the  coefficient ring
$E(n)_\ast= \mathbb{Z}_{(p)}[v_1,v_2,\ldots,v_{n-1}][v_n^{\pm
1}]$.
 The $MU_\ast$-algebra structure on $E(n)_\ast$
 is defined by the correspondence $x_i\rightarrow 0$ for $i\neq 2(p^n-1)$ and
$x_{2(p^n-1)}\rightarrow v_n$ where $x_i$ are the Hazewinkel
generators of $MU_{(p)\ast}$. We have the following
\begin{cor}\label{yu}
There exist  $MU$-algebras $\hat{E}(n)$, $\widehat{EL}(n)$,
$K(n)$, $KL(n)$, and the homotopy commutative diagram of
$MU$-algebras
\begin{equation}\label{gt}\xymatrix{
\hat{E}(n)\ar[d]\ar[r]&\widehat{EL}(n)\ar[d]\\
K(n)\ar[r]&KL(n)}\end{equation} which realizes in homotopy the
diagram of $MU_\ast$-algebras
$$\xymatrix{\hat{\mathbb{Z}_p}[v_n^{\pm
1}][[v_1,v_2,\ldots,v_{n-1}]]\ar[d]\ar[r]& W(L)[v_n^{\pm
1}][[v_1,v_2,\ldots,v_{n-1}]]\ar[d]\\ \mathbb{F}[v_n^{\pm
1}]\ar[r]& L[v_n^{\pm 1}]}$$
\end{cor}
{\bf Proof.} We first construct the map of $MU_{(p)}$-algebras
$k(n)\rightarrow kL(n)$ by adjoining $v_n$ to the extension
$H\mathbb{F}_p\rightarrow HL$. Here $k(n)$ and $kL(n)$ are the
connective Morava $K$-theories with coefficient rings $k(n)_\ast =
\mathbb{F}_p[v_n]$ and $kL(n)_\ast = L[v_n]$. Adjoining
indeterminates $p,v_1,v_2,\ldots, v_{n-1}$ we construct the
homotopy commutative diagram of $MU$-algebras
\begin{equation}\label{df}\xymatrix{{\widehat{BP}\langle
n\rangle}\ar[d]\ar[r]& \widehat{BPL}\langle
n\rangle\ar[d]\\k(n)\ar[r]& kL(n)}\end{equation} which realizes in
homotopy the diagram
$$\xymatrix{\hat{\mathbb{Z}_p}[v_n][v_1,v_2,\ldots,v_{n-1}]\ar[d]\ar[r]&
W(L)[v_n][v_1,v_2,\ldots,v_{n-1}]\ar[d]\\
\mathbb{F}_p\ar[r]&L[v_n]}$$

According to Proposition \ref{oj} Bousfield localization of
${\widehat{BP}\langle n\rangle}$   with respect to $K(n)$ on the
level of coefficient rings amounts to inverting $v_n$ and
completing at the ideal $(v_1,v_2,\ldots, v_{n-1})$. Further
notice that since $\pi_\ast {\widehat{BPL}\langle n\rangle}$ is a
free $\pi_\ast {\widehat{BP}\langle n\rangle}$-module of finite
rank (which is equal to the degree of the field extension
$L/\mathbb{F}_p$) the $MU$-algebra ${\widehat{BPL}\langle
n\rangle}$ is  a finite cell $\widehat{{BP}}\langle
n\rangle$-module and therefore its $K(n)$-localization is
equivalent to ${\widehat{BP}\langle
n\rangle}_{K(n)}\wedge_{\widehat{BP}\langle
n\rangle}{\widehat{BPL}\langle n\rangle}$ and its coefficient ring
is $$W(L)\otimes\pi_\ast{\widehat{BP}\langle n\rangle}_{K(n)}=
W(L)[v_n^{\pm 1}][[v_1,v_2,\ldots,v_{n-1}]].$$ Therefore we obtain
the diagram (\ref{gt}) by localizing (\ref{df}) with respect to
$K(n)$ in the category of $MU$-algebras. Corollary \ref{yu} is
proved.

\section{Computing homotopy groups of mapping spaces}
In this section we investigate spaces of strictly multiplicative
maps from $\hat{EL}(n)$ to certain $MU$-algebras which we call
called `strongly $KL(n)$-complete'. As a consequence we obtain
versions of the Hopkins-Miller theorem as well as splitting
theorems for such spectra. Such theorems (in a weaker, up to
homotopy form) were previously obtained by methods of formal group
theory, cf. \cite{ BaW}. Recall that $\widehat{EL}(n)$ and $KL(n)$
denote topological extensions of spectra $\hat{E}(n)$ and $K(n)$
corresponding to the separable field extension
$\mathbb{F}_p\hookrightarrow L$. All our results remain valid and
are still nontrivial for $L=\mathbb{F}_p$ in which case
$\widehat{EL}(n)$ and $KL(n)$ specialize to the completed
Johnson-Wilson theory $\hat{E}(n)$ and Morava $K$-theory $K(n)$
respectively.

Let $A$, $B$ be $S$-algebras, which we will assume without loss of
generality to be $q$-cofibrant in the sense of \cite{EKMM}. That
means in particular, that the topological space of $S$-algebra
maps $B\rightarrow A$ has the `correct' homotopy type (i.e. the
one that depends only on the homotopy type of $A$ and $B$ as
$S$-algebras). Denote this topological space by $F_{S-alg}(B,A)$.
Also denote the set of multiplicative up to homotopy maps from $A$
to $B$ by $Mult (A,B)$.

We recall one result from \cite{Laz} which will be needed later on.
\begin{theorem}\label{ext} Let $\Sigma^{-1}M\rightarrow X\rightarrow A$ be a singular
extension of $S$-algebras associated with a derivation
$d:A\rightarrow A\vee M$ and $f:B\rightarrow A$ a map of
$S$-algebras. Then $f$ lifts to an $S$-algebra map $B\rightarrow
X$ iff a certain element in $Der^0(B,M)$ is zero. Assuming that a
lifting exists  the homotopy fibre of the map
\begin{equation}\label{po}F_{S-alg}(B,X)\rightarrow
F_{S-alg}(B,A)\end{equation} over the point $f\in F_{S-alg}(B,A)$
is weakly equivalent to $\Omega^\infty {\bf Der}(B,\Sigma^{-1}M)$,
the $0$th space of the spectrum ${\bf Der}(X,\Sigma^{-1}M)$. In
particular if the spectrum ${\bf Der}(X,\Sigma^{-1}M)$ is
contractible then (\ref{po}) is a weak equivalence.
\end{theorem}

Recall that the spectrum  $\widehat{EL}(n)$ with coefficient rings
$W(L)[v_n^{\pm 1}][[v_1,v_2,\ldots,v_{n-1}]]$ has a structure of
an $S$-algebra.

We say that the generalized Hopkins-Miller theorem holds for an
$S$-algebra $E$ if the following two conditions
hold:\begin{enumerate}\item $\pi_0F_{S-alg}(\widehat{EL}(n),
E)=Mult(\widehat{EL}(n), E)$ \item
$\pi_iF_{S-alg}(\widehat{EL}(n), E)=0$ for $i>0$.\end{enumerate}

\begin{prop}\label{bk} The generalized Hopkins-Miller theorem holds for the
$S$-algebra $KL(n)$.
\end{prop}
{\bf Proof.} Consider the Bousfield-Kan spectral sequence
(cf.\cite{Bou}) for the mapping space $F_{S-alg}(\widehat{EL},E)$.
The identification of the $E_2$-term  is standard and we refer the
reader to \cite{Rez} for necessary details. Here we only mention
that the key ingredient in this identification is the existence of
the K\"unneth formula $$KL(n)_*(\widehat{E}(n)\wedge
\widehat{EL}(n))= KL(n)_*\widehat{EL}(n)
\otimes_{KL(n)_*}KL(n)_*\widehat{EL}(n).$$ This is the result we
need:$$E_2^{st}=
\begin{cases}
Der^{st}_{KL(n)_*}(\widehat{EL}(n)_*K(n), KL(n)_*)&\text {for
$(s,t)\neq (0,0)$},
\\ Mult(\widehat{EL}(n), KL(n) )& \text{for $(s,t)=(0,0)$}.
\end{cases}$$
Here $Der_k^{\ast\ast}(?,?)$ is defined for a graded $k$-algebra
$R_\ast$ and a graded
 $R_\ast$-bimodule $M_\ast$ as the shifted Hochschild cohomology:
$$Der_k^{\ast\ast}(A_\ast,M_\ast):=HH_k^{\ast+1\ast}(A_\ast,M_\ast)$$
(the second grading reflects the fact that $A_\ast$ and $M_\ast$ are graded
objects).

Further $$\widehat{EL}(n)_* KL(n)=\hat{E}(n)_\ast K(n)\otimes
W(L)\otimes L= E(n)_\ast K(n)\otimes L\otimes L$$ and computations
in \cite{Rav}, Chapter VI show that $$E(n)_\ast
K(n)=\Sigma_\ast=\mathbb{F}_p[v_n^{\pm
1}][t_k|k>0]/(t_k^{p^n}-v_n^{p^k-1}t_k)$$ where the degree of
$t_k$ is $2(p^k-1)$. We have:
$$HH^{\ast\ast}_{KL(n)_*}(\Sigma_\ast\otimes L\otimes L,KL(n)_*)=
HH^{\ast\ast}_{KL(n)_*}(\Sigma_*\otimes L,K(n)_\ast\otimes L)$$
$$=HH^{\ast\ast}_{K(n)_*}(\Sigma_*, K(n)_\ast)\otimes
HH^\ast_{\mathbb{F}_p}(L,L)=
HH^{\ast\ast}_{K(n)_\ast}(\Sigma_\ast, K(n)_\ast)\otimes L.$$
 An easy cohomological calculation
(due to A.Robinson, cf. \cite{Rob}) shows that
$HH^{\ast\ast}_{K(n)_\ast}(\Sigma_\ast, K(n)_\ast)=K(n)_\ast$ and therefore
\begin{equation}\label{bnm}
HH^{\ast\ast}_{KL(n)_*}(\widehat{EL}(n)_* KL(n))=K(n)\otimes L=
KL(n)_*\end{equation} Thus the Bousfield-Kan spectral sequence
reduces to its corner term and Proposition \ref{bk} is proved.
\begin{rem} Notice that even though in Proposition \ref{bk} we used
the results of the previous section that spectra $KL(n)$,  and
$\widehat{EL}(n)$ admit $S$-algebra structures we did not specify
which structures are used. In other words Proposition \ref{bk} is
valid with {\it any} choice of $S$-algebra structures on the
spectra in question.
\end{rem}
\
We can now regard $KL(n)$  as a bimodule over $\widehat{EL}(n)$ by
choosing an $S$-algebra map $\widehat{EL}(n)\rightarrow KL(n)$
(which is
 supplied by
Proposition \ref{bk}). Therefore it makes sense to consider
spectra of topological Hochschild cohomology and topological
derivations of $\widehat{EL}(n)$ with values in $KL(n)$. It turns
out that these spectra do not depend up to weak equivalence which
bimodule structure we choose.
\begin{prop} \label{spl} The canonical map of $S$-modules
$${\bf THH}(\widehat{EL}(n),KL(n))\rightarrow KL(n)$$ is a  weak
equivalence. Furthermore the $S$-module ${\bf
Der}(\widehat{EL}(n),KL(n))$ is contractible.
\end{prop}
{\bf Proof.} Since there is a homotopy cofibre sequence of
$S$-modules
 $$\Sigma^{-1}{\bf Der}(\widehat{EL}(n),KL(n))\rightarrow
{\bf THH}(\widehat{EL}(n),KL(n))\rightarrow KL(n)$$ it suffices to
prove the statement about $\bf THH$. Consider the following
spectral sequence
\begin{equation}\label{ful}
Ext^{\ast\ast}_{\widehat{EL}(n)_*\widehat{EL}(n) }
(\widehat{EL}(n)_* , KL(n)_*) \Rightarrow
 THH^\ast(\widehat{EL}(n),KL(n))\end{equation}
We have the following isomorphisms:
$$Ext^{\ast\ast}_{\widehat{EL}(n)_{\ast}\widehat{EL}(n) }
(\widehat{EL}(n)_{\ast} , KL(n)_{\ast})$$ $$\cong
Ext^{\ast\ast}_{\widehat{EL}(n)_{\ast}\widehat{EL}(n)\otimes_{\widehat{EL}(n)_*}KL(n)_*
} (\widehat{EL}(n)_{\ast}\otimes_{\widehat{EL}_*}KL(n)_* ,
KL(n)_{\ast})$$ $$ \cong Ext^{\ast\ast}_{\widehat{EL}(n)_{\ast}
KL(n) } (KL(n)_* , KL(n)_{\ast})\cong
HH^\ast_{KL(n)_*}(\widehat{EL}(n)_* KL(n), KL(n)_*).$$ The first
isomorphism follows from the fact that $E(n)_*$-module
$\widehat{EL}(n)_*\widehat{EL}(n)$ is flat (\cite{Rez},
Proposition 15.1) and the second and third isomorphisms are
obvious.

 The isomorphism (\ref{bnm})   shows that the
spectral sequence (\ref{ful}) collapses giving the weak
equivalence $${\bf THH}(\widehat{EL}(n),KL(n))\simeq KL(n)$$
  and
Proposition \ref{spl} is proved.

We now describe a particular class of $MU$-algebras $E$ which will
be called strongly $KL(n)$-complete and for which the set of
multiplicative maps~ $\widehat{EL}(n)\rightarrow E$ has an
especially simple form.

\begin{defi}
The  category $SC(KL(n))$ of strongly $KL(n)$-complete
$MU$-algebras is the smallest subcategory of the homotopy category
of $MU$-algebras over $KL(n)$ which contains $KL(n)$ itself and
closed under singular extensions with fibre $KL(n)$ and inverse
limits of towers of $MU$-algebras over $KL(n)$.
\end{defi}
\begin{rem} Note that the $MU_*$-module $\pi_*KL(n)$ is linearly
compact and the tower of linearly compact modules has vanishing
$lim^1$. Since linearly compact modules are closed under inverse
limits  and singular extensions we conclude that for any $E\in
SC(KL(n))$ the $MU_*$-module $E_*$ is linearly compact. In
particular we never have to worry about $lim^1$-problems.
Furthermore it follows that the coefficient rings of strongly
$KL(n)$-complete $MU$-algebras are concentrated in even degrees.
\end{rem}\begin{rem} Examples of strongly $KL(n)$ - complete $MU$-algebras
include $KL(n)$ and $\widehat{EL}(n)$. To extend the number of
examples recall the notion of the {\it Artinian completion} of a
spectrum due to A.Baker and U.W\"urgler, cf. \cite{BaW}. Let $R$
be a commutative $S$-algebra and $I$ be a graded maximal ideal in
the coefficient ring of $R$ generated by a regular sequence
$u_1,u_2,\ldots$, possibly infinite. Then for an $R$-module $M$
define the Artinian completion $\hat{M}$ of $M$ as the homotopy
inverse limit $$\hat{M}:=\holim
M\wedge_RR/(u_n^{i_n},u_{n+1}^{i_{n+1}},\ldots)$$ where the
inverse limit is taken over the collections of positive integers
$n,i_n,i_{n+1},\ldots$. In fact Baker and W\"urgler have a more
general definition of the Artinian completion but the one just
given is sufficient for our purposes.

Now taking $R=MU_{(p)}$ and
$M=v_n^{-1}BPL:=v_n^{-1}[HL(v_1^\infty,v_2^{\infty}\ldots)]$ or
$M=v_n^{-1}PL(n):=v_n^{-1}BPL/(v_1,v_2,\ldots,v_{n-1})$ we could
form Artinian completions $\widehat{v_n^{-1}BPL}$ and
$\widehat{v_n^{-1}PL(n)}$. It follows easily that
$\widehat{v_n^{-1}BPL}$ and $\widehat{v_n^{-1}PL(n)}$ are strongly
$KL(n)$-complete $MU$-algebras.\end{rem}

Let $E$ be a strongly $KL(n)$-complete  $MU$-algebra. The
canonical $S$-algebra map (which is also an $MU$-algebra map)
$E\rightarrow KL(n)$
 determines via composition the
morphism of topological spaces (unbased)
\begin{equation}\label{sap}
F_{S-alg}({EL}(n),E)\rightarrow F_{S-alg} ({EL}(n),KL(n))
\end{equation}
\begin{prop}\label{vb} The map (\ref{sap}) is a weak equivalence.
In particular $F_{S-alg}({EL}(n),E)$ is a homotopically discrete
space. Moreover the $S$-module ${\bf Der}({EL}(n),E)$ is
contractible.
\end{prop}
{\bf Proof.} We filter the class $SC(KL(n))$ by subclasses
$\{SC_\nu(KL(n))\}$ where $\nu$ is an ordinal. The lowest level
$SC_0(KL(n))$ consists of one copy of $KL(n)$. Suppose the class
$\{SC_\nu(KL(n))\}$ has been defined for $\nu<\mu$. We say that an
$MU$-algebra $E$ belongs to the class $SC_{\mu}(KL(n))$ if and
only if
\begin{enumerate}\item $E$  weakly equivalent to an inverse limit of
the tower of $MU$-algebras over $KL(n)$ with each term belonging
to some $SC_{\nu}(KL(n))$, $\nu<\mu$ or
\item there is a singular extension
$KL(n)\rightarrow E\rightarrow F$ where $F\in SC_{\nu}(KL(n))$ and
$\nu<\mu$.
\end{enumerate} To prove that the map (\ref{sap}) is a weak equivalence we
will use transfinite induction up the filtration on $SC(KL(n))$.
The base case $\nu=0$ is implied by Propositions \ref{bk} and
\ref{spl}. Suppose that our statement is true for $MU$-algebras
belonging to $SC_\nu(KL(n))$ for $\nu<\mu$ and consider $E\in
SC_{\mu}(KL(n))$. If $E$ can be included in the singular extension
as in the condition (2) above then applying Theorem \ref{ext} and
Proposition \ref{spl} we get the desired conclusion.

Finally suppose that $E$ is the inverse limit of the tower
$\{E_l\}$ where $E_l\in SC_{\nu_l}(KL(n))$ and $\nu_l<\mu$.
Applying the functor $F_{S-alg}(\widehat{EL}(n),?)$ to $\{E_l\}$
we get a tower of topological spaces
$\{F_{S-alg}(\widehat{EL}(n),E_l)\}$. By inductive assumption all
spaces $F_{S-alg}(\widehat{EL}(n),E_l)$ are weakly equivalent to
$F_{S-alg}(E,KL(n))$ and the maps in the tower
$\{F_{S-alg}(\widehat{EL}(n),E_l)\}$ respect this equivalence. We
conclude that the space $F_{S-alg}(\widehat{EL}(n),E)$ is also
weakly equivalent to $F_{S-alg}(\widehat{EL}(n),KL(n))$. Similar
arguments show that ${\bf Der}(\widehat{EL}(n),KL(n))$  is
contractible and Proposition \ref{vb} is proved.

\begin{rem} We see that for any strongly $KL(n)$-complete $MU$-algebra $E$
the canonical $S$-algebra map $\widehat{EL}(n)\rightarrow KL(n)$
lifts uniquely to a map $\widehat{EL}(n)\rightarrow E$. Therefore
$E$ is naturally a (bi)module over $\widehat{EL}(n)$.
\end{rem}
The following lemma descibes the structure of cohomology
operations from $\widehat{EL}(n)$ to a strongly $KL(n)$-complete
$MU$-algebra $E$.
\begin{lem}\label{stru} Let $E$ be a  strongly $KL(n)$-complete $MU$-algebra. Then
the evaluation map $$E^*\widehat{EL}(n)\rightarrow
Hom_{\widehat{EL}(n)_*}(\widehat{EL}(n)_* \widehat{EL}(n), E_*)$$
is an isomorphism.
\end{lem}
{\bf Proof.} We have the following natural isomorphism of
$S$-modules $$F_S(\widehat{EL}(n),E)\cong F_{\widehat{EL}(n)}
(\widehat{EL}(n)\wedge \widehat{EL}(n) ,E)$$ Consider further the
spectral sequence
\begin{eqnarray}\label{spe} \nonumber Ext^{\ast\ast}_{\widehat{EL}(n)_{\ast}}
(\widehat{EL}(n)_{\ast} \widehat{EL}(n),E_\ast) \Rightarrow
[\widehat{EL}(n)\wedge \widehat{EL}(n) ,E]^\ast_{\widehat{EL}(n)}=
[\widehat{EL}(n),E]^\ast\end{eqnarray} We claim that all higher
$Ext$ groups vanish so that our spectral sequence reduces to its
first column which is $Ext^{0\ast}_{\widehat{EL}(n)_{\ast}}
(\widehat{EL}(n)_{\ast} \widehat{EL}(n), E_\ast)$. (This would
clearly give us the statement of the lemma). To see this first
assume that $E=KL(n)$

 Since the
$\widehat{EL}(n)_*$-module $\widehat{EL}(n)_{\ast}
\widehat{EL}(n)$ is flat we have by flat base change:
$$Ext^{i\ast}_{\widehat{EL}(n)_{\ast}}
(\widehat{EL}(n)_*\widehat{EL}(n),KL(n)_*)=Ext^{i\ast}_{KL(n)_*}
(\widehat{EL}(n)_*
\widehat{EL}(n)\otimes_{\widehat{EL}(n)_*}KL(n)_*,KL(n)_*)$$ and
the last $Ext$-group is zero for $i>0$ since $KL(n)_*$ is
injective as a graded module over itself.

Therefore our claim about the vanishing of  higher $Ext$-groups in
(\ref{spe}) is proved for  $E=KL(n)$. To get the general case we
will use transfinite induction as in Proposition \ref{vb}.

Suppose that $E$ can be included in a singular extension
\begin{equation}\label{bhu}KL(n)\rightarrow E\rightarrow F\end{equation}
 where $F\in SC(KL(n))$ and
such that
$Ext^{i\ast}_{\widehat{EL}(n)_*}(\widehat{EL}(n)_*\widehat{EL}(n)_*,F_*)=0$
for $i>0$. Since the coefficient rings of spectra $KL(n), E,F$ are
even the homotopy long exact sequence corresponding to the cofibre
sequnce (\ref{bhu}) is in fact a {\it short} exact sequence of
coefficient rings:
\begin{equation}\label{bhu1}KL(n)_*\rightarrow E_*\rightarrow
F_*\end{equation} Associated to (\ref{bhu1}) is a long exact
sequence
$$Ext^{1\ast}_{\widehat{EL}(n)_*}(\widehat{EL}(n)_*\widehat{EL}(n)_*,KL(n)_*)\rightarrow
Ext^{1\ast}_{\widehat{EL}(n)_*}(\widehat{EL}(n)_*\widehat{EL}(n)_*,E_*)$$
$$ \rightarrow
Ext^{1\ast}_{\widehat{EL}(n)_*}(\widehat{EL}(n)_*\widehat{EL}(n)_*,F_*)\rightarrow
Ext^{2\ast}_{\widehat{EL}(n)_*}(\widehat{EL}(n)_*\widehat{EL}(n)_*,KL(n)_*)\rightarrow\ldots$$
from which it follows that
$Ext^{i\ast}_{\widehat{EL}(n)_*}(\widehat{EL}(n)_*\widehat{EL}(n)_*,E_*)=0$.

 Further let $E$ be the homotopy inverse
limit of strongly $KL(n)$-complete $MU$-algebras $E_l$ for which
higher $Ext$'s do vanish. Then $$Ext^{i\ast}_{\widehat{EL}(n)_*}
(\widehat{EL}(n)_* \widehat{EL}(n),E_*)
=Ext^{i\ast}_{\widehat{EL}(n)_*} (\widehat{EL}(n)_*
\widehat{EL}(n),\lim_{l\rightarrow \infty}E_{l\ast})$$
$$=\lim_{l\rightarrow\infty}Ext^{i\ast}_{\widehat{EL}(n)_*}
(\widehat{EL}(n)_* \widehat{EL}(n),E_{l\ast})=0$$ for $i>0$. This
concludes the proof of Lemma \ref{stru}.

We can now formulate our main
\begin{theorem} \label{ho}The generalized Hopkins-Miller theorem holds for any
strongly $KL(n)$-complete $MU$-algebra $E$.
\end{theorem}
{\bf Proof.} By Proposition \ref{vb} we know that the space
$F_{S-alg}(\widehat{EL}(n),E)$ is homotopically discrete. We need
to check that the map
\begin{equation}\label{zab}\pi_0F_{S-alg}(\widehat{EL}(n),E)\rightarrow
Mult(\widehat{EL}(n),E)\end{equation} is a bijection. Consider the
following commutative diagram
$$\xymatrix{\pi_0F_{S-alg}(\widehat{EL}(n),E)\ar[d]\ar[r]&
\pi_0F_{S-alg}(\widehat{EL}(n),KL(n))\ar[d]\\Mult(\widehat{EL}(n),E)\ar[r]&
Mult(\widehat{EL}(n),KL(n))}$$
 By Propostions \ref{bk}
and \ref{vb} the upper horizontal and the right vertical arrows
are bijections. It follows that the map (\ref{zab}) is injective.
So it remains to show that an arbitrary multiplicative operation
$\widehat{EL}(n)\rightarrow E$ lifts to an $S$-algebra map. Since
any such operation determines an $\widehat{EL}(n)_{\ast}$-algebra
map $\widehat{EL}(n)_{\ast} \widehat{EL}(n)\rightarrow E_\ast$ we
conclude by Lemma \ref{stru} that there is an injective map
\begin{equation}\label{dfg} Mult(\widehat{EL}(n),E)\rightarrow
Hom_{\widehat{EL}(n)_{\ast}-alg} (\widehat{EL}(n)_{\ast}
\widehat{EL}(n),E_\ast)\end{equation} We claim that any
$\widehat{EL}(n)_{\ast}$-algebra map $\widehat{EL}(n)_{\ast}
\widehat{EL}(n)\rightarrow KL(n)_{\ast}$ lifts uniquely to a map
$\widehat{EL}(n)_{\ast} \widehat{EL}(n)\rightarrow E_{\ast}$ so
that there is an isomorphism $$\label{rte}
Hom_{\widehat{EL}(n)_{\ast}-alg}(\widehat{EL}(n)_{\ast}\widehat{EL}(n),E_\ast)
\cong Hom_{\widehat{EL}(n)_{\ast}-alg}(\widehat{EL}(n)_{\ast}
\widehat{EL}(n),KL(n)_{\ast}).$$ Using inductive arguments as in
the previous lemma we see that  our claim is equivalent to the
vanishing of certain Hochschild cohomology classes in $HH^{\ast 2}
_{\widehat{EL}(n)_{\ast}}(\widehat{EL}(n)_{\ast}
\widehat{EL}(n),KL(n)_{\ast})$. Flat base change gives an
isomorphism $$HH^{\ast\ast}_{\widehat{EL}(n)_{\ast}}
(\hat{EL}(n)_{\ast} \widehat{EL}(n),KL(n)_{\ast})$$
$$=HH^{\ast\ast}_{KL(n)_{\ast}} (\widehat{EL}(n)_{\ast}
\widehat{EL}(n)\otimes_{\widehat{EL}(n)_{\ast}}KL(n)_{\ast}
,KL(n)_{\ast})$$
$$=HH^{\ast\ast}_{KL(n)_{\ast}}(\widehat{EL}(n)_*KL(n)
,KL(n)_{\ast})= KL(n)_{\ast}.$$ We see that all higher Hochschild
cohomology groups are zero and therefore our claim is proved.

Further we have an obvious change of rings isomorphism
$$Hom_{\widehat{EL}(n)_{\ast}-alg} (\widehat{EL}(n)_{\ast}
\widehat{EL}(n),KL(n)_{\ast}) \cong
Hom_{KL(n)_{\ast}-alg}(\widehat{EL}(n)_*KL(n),KL(n)_{\ast})$$ The
last term corresponds bijectively by Proposition \ref{vb} to
homotopy classes of $S$-algebra maps from $\widehat{EL}(n)$ to
$E$. This shows that any multiplicative operation
$\widehat{EL}(n)\rightarrow E$ lifts to an $S$-algebra map and
Theorem \ref{ho} is proved.
\begin{rem} The class of strongly $KL(n)$-complete $MU$-algebras is probably
not the most general one for which the conclusion of Theorem
\ref{ho} holds. It is plausible that one has the generalized
Hopkins-Miller theorem for any $KL(n)$-local $MU$-algebra but such
a result seems out of reach at present.
\end{rem}
 Tracing the proofs of Proposition \ref{bk} and Theorem \ref{ho}
we see that they don't depend on the choice of the $S$-algebra
structure on $\widehat{EL}(n)$ as long as the ring spectrum
structure is fixed. (But they do depend on the choice of the
$S$-algebra structure on $E$.) This gives the following
\begin{cor}
The spectrum $\widehat{EL}(n)$ has a unique (up to a noncanonical
isomorphism)
 structure of an $S$-algebra compatible
with its structure of a ring spectrum.
\end{cor}
{\bf Proof.} Let $\widehat{EL}^\prime(n)$ denote the $S$-algebra
whose underlying ring spectrum is equivalent to $\widehat{EL}(n)$
but the $S$-algebra structure is possibly different. Then the
given multiplicative up to homotopy weak equivalence
$\widehat{EL}^\prime(n)\rightarrow \widehat{EL}(n)$ can be lifted
to an $S$-algebra map which shows that the would-be exotic
$S$-algebra structure on $\widehat{EL}^\prime(n)$ is actually
isomorphic to the standard one.

\begin{rem} A similar statement about the uniqueness of the $A_\infty$ structure
on $\hat{E}(n)$ was proved by A.Baker in \cite{Bak}. However the proof in the
cited reference was based on the obstruction theory of A.Robinson and it is unclear
at this time  to what extent this theory carries over in the present context.
\end{rem}
Theorem \ref{ho} allows one to obtain splittings of various
strongly $KL(n)$-complete $MU$-algebras. For example consider
$\widehat{v_n^{-1}BPL}$, the Artinian completion of the
$v_n$-localization of the Brown-Peterson spectrum $BPL$. Then
$\widehat{v_n^{-1}BPL}$ is a strongly $KL(n)$-complete
$MU$-algebra obtained from $KL(n)$ by adjoining the collection of
indeterminates corresponding to the regular sequence
$\{p,v_1,v_2,\ldots,v_{n-1},v_{n+1},\ldots\}$. There is a
canonical map $\widehat{v_n^{-1}BPL}\rightarrow \widehat{EL}(n)$
which on the level of coefficient ring reduces to killing the
ideal $(v_{n+1},v_{n+2},\ldots)$ in $\widehat{v_n^{-1}BPL}_*$.
From Theorem \ref{ho} we conclude that any $S$-algebra map
$\widehat{EL}(n)\rightarrow \widehat{EL}(n)$ lifts uniquely to an
$S$-algebra map $\widehat{EL}(n)\rightarrow
\widehat{v_n^{-1}BPL}$. In particular the identity map
$id:\widehat{EL}(n)\rightarrow \widehat{EL}(n)$ can be so lifted
and we obtain
\begin{theorem}\label{suk} The $S$-algebra map $\widehat{v_n^{-1}BPL}\rightarrow
\widehat{EL}(n)$ admits a unique $S$-algebra splitting
$\widehat{EL}(n)\rightarrow \widehat{v_n^{-1}BPL}$
\end{theorem}
An analogous theorem was proved (in a weaker, up to homotopy form)
in \cite{BaW}.

\end{document}